\documentclass[11pt]{article}
\usepackage{amsfonts}
\usepackage{amssymb}

\newtheorem{theorem}{Theorem}[section]
\newtheorem{proposition}[theorem]{Proposition}
\newtheorem{lemma}[theorem]{Lemma}
\newtheorem{corollary}[theorem]{Corollary}

\newtheorem{remark}{Remark}[section]
\newtheorem{example}{Example}[section]
 \usepackage{amssymb}
\usepackage{epsfig}
\usepackage{amsmath}
\usepackage{amsfonts}

\newcommand{\R}{{{\Bbb R}}}

\newcommand{\N}{{{\Bbb N}}}

\def\qed{\hbox to 0pt{}\hfill$\rlap{$\sqcap$}\sqcup$\medbreak}

\title{Mean value integral inequalities\footnote{Partially
supported by FEDER and
Ministerio de Educaci\'on y Ciencia, Spain,
project MTM2010-15314.}}
\date{June, 2011}
\begin{document}
\maketitle

\vspace{1cm}

\begin{center}
{\large Rodrigo L\'opez Pouso \\
Department of Mathematical Analysis\\
Faculty of Mathematics,\\University of Santiago de Compostela,\\
15782 Santiago de
Compostela, Spain.
%\\
%Phone: 34 981 56 31 00,
%Ext. 13213 \\ FAX: 34 981 59 70 54
}
\end{center}

\begin{abstract}
Let $F:[a,b]\longrightarrow \R$ have zero derivative in a dense subset of $[a,b]$. What else we need to conclude that $F$ is constant in $[a,b]$? We prove a result in this direction using some new Mean Value Theorems for integrals which are the real core of this paper. These Mean Value Theorems are proven easily and concisely using Lebesgue integration, but we also provide alternative and elementary proofs to some of them which keep inside the scope of the Riemann integral.
\end{abstract}

 \section{Introduction}
Roughly, zero derivative implies constancy, but the devil hides in the details. This note follows the spirit of \cite{koliha}, where we find the following theorem:
\begin{theorem}
\label{tkol}
A function $F:I=[a,b] \longrightarrow \R$ is constant in $I$ provided that one of the following conditions holds:
\begin{enumerate}
\item[(A)] $F'=0$ everywhere in $I$; or
\item[(B)] $F'=0$ nearly everywhere in $I$ (i.e., $F'=0$ in $[a,b]\setminus C$, with $C$ countable) and $F$ is continuous; or
\item[(C)]$F'=0$ almost everywhere in $I$ and $F$ is absolutely continuous.
\end{enumerate}
\end{theorem}

We observe in conditions (A), (B) and (C) that the bigger the exceptional set is, the more regular the function $F$ must be. Notice that any of the conditions (A), (B) and (C) implies that $F'=0$ in a dense subset of $[a,b]$, which makes us wonder how $F$ should be in that case in order to ensure that $F$ is constant. In this paper we prove the following result.
\begin{theorem}
\label{t1}
A function $F:I=[a,b]\longrightarrow \R$ is constant in $I$ provided that the following condition holds:
\begin{enumerate}
\item[(D)] $F'=0$ densely in $I$ (i.e., $F'$ exists and is equal to zero in a dense subset of $I$) and $F$ is an indefinite Riemann integral, i.e., we can find a Riemann--integrable function $f:I\longrightarrow \R$ and a constant $c \in \R$ such that
$$F(x)=c+\int_a^x{f(y) \, dy} \quad (x \in I).$$
\end{enumerate}
\end{theorem}  

We follow the spirit of Koliha's paper \cite{koliha} in the sense that we are going to present some very mean Mean Value Theorems (according to the terminology in \cite{koliha}, the meaner, the stronger), and then we will use them to prove part (C) in Theorem \ref{tkol} and Theorem \ref{t1}. In fact, these Mean Value Theorems are the main subject in our paper, and they have some other interesting consequences. In particular, we deduce a very mean Mean Value Theorem for the Riemann integral which leads to a short proof for our Theorem \ref{t1} in section 4. 

\section{A first proof of Theorem \ref{t1} and some remarks}

Our first proof is based on the fact that (D) implies (C), and then Theorem \ref{t1} is a particular case to part (C) in Theorem \ref{tkol}. 
\bigbreak

\noindent
{\bf First proof of Theorem \ref{t1}.} Let $F$ and $f$ be as in the statement. We start proving that if $f$ is continuous at $x \in (a,b)$ then $f(x)=0$. Reasoning by contradiction, assume that $f$ is continous at $x \in (a,b)$ and $f(x) \neq 0$. This ensures the existence of some $r>0$ such that $|f|>r$ everywhere in $(x-r,x+r)$ and $f$ has constant sign in $(x-r,x+r)$. Let $z \in (x-r,x+r)$ be such that $F'(z)=0$. For all $y \in (z,x+r)$, $y$ sufficiently close to $z$, we have
$$r>\left|\dfrac{F(y)-F(z)}{y-z}\right|=\dfrac{1}{y-z}\left|\int_z^y{f(s) \, ds}\right| \ge r,$$
a contradiction. Therefore $f(x)=0$ whenever $f$ is continuous at $x \in (a,b)$, and then the Fundamental Theorem of Calculus and Lebesgue's test guarantee that $F'=f=0$ almost everywhere in $I$. The conclusion now follows from part (C) in Theorem \ref{tkol}. \qed

\noindent
{\bf Remarks to Theorem \ref{t1} and its first proof:}
\begin{enumerate}
\item Indefinite Riemann integrals are for Riemann integration what absolutely continous functions are for the Lebesgue integral. 

Indefinite Riemann integrals have been recently characterized in \cite{tho} as follows: a function $F:[a,b]\longrightarrow \R$ is an indefinite Riemann integral if and only if for all $\varepsilon >0$ a positive $\delta$ can be found so that
$$\sum_{i=1}^n\left| \dfrac{F(\xi_i)-F(x_{i-1})}{\xi_i-x_{i-1}}- \dfrac{F(x_i)-F(\xi_i')}{x_i-\xi_{i}'}\right|
(x_i-x_{i-1})<\varepsilon$$
for every subdivision $a=x_0<x_1<\cdots<x_n=b$ that is finer than $\delta$ and every choice of associated points $x_{i-1}<\xi_i \le \xi_i'<x_{i}$.

\item Indefinite Riemann integrals are Lispchitz, but the converse is false and Theorem \ref{t1} is not valid with indefinite Riemann integrals replaced by Lispchitz functions. The following example justifies these two statements. 

\begin{example}
Let $C \subset [0,1]$ be a Cantor set with Lebesgue measure $\mu >0$ and let $\chi_C$ be its characteristic function. Using Lebesgue integral we define 
$$F(x)=\int_0^x{\chi_C(s) \, ds} \quad (x \in [0,1]).$$
Notice that $F$ is Lispchitz on $[0,1]$, $F'=0$ in $[0,1]\setminus C$ (hence $F'=0$ densely in $[a,b]$), and
$$F(1)=\int_0^1{\chi_C(s) \, ds}=\mu >0=F(0).$$
\end{example}
Lipschitz functions are absolutely continuous, so Example 1.1 also shows that Theorem \ref{t1} is not valid with indefinite Riemann integrals replaced by indefinite Lebesgue integrals (i.e., absolutely continuous functions).

\item The proof of Theorem \ref{t1} reveals that condition (D) implies condition (C), so (D) is not an essentially new situation. This raises the problem of determining a set of functions for which the condition $F'=0$ densely in $[a,b]$ implies that $F$ is constant but does not imply that $F$ satisfies (C). We know from Theorem \ref{t1} and the previous remark that such set (if it exists!) should be bigger than the set of indefinite Riemann integrals and smaller than the set of Lipschitz functions.
\end{enumerate}

To close this section let us point out that our first proof of Theorem \ref{t1} leans on null--measure sets and absolutely continuous functions, which neither are present in the statement nor seem naturally connected with the assumptions. Can we have another proof which does not use these elements? The answer is positive and we give one such proof in Section 4.

\section{Mean Value Theorems for Lebesgue integrals}
Part (C) in Theorem \ref{tkol} can be proven with the aid of the following Mean Value Theorem. In the sequel $m$ stands for the Lebesgue measure in $\R$.
\begin{theorem}
\label{vmmvt}
If $f:I=[a,b]\longrightarrow \overline{\R}$ is Lebesgue--integrable then
\begin{equation}
\label{pm1}
m\left(\left\{c \in (a,b) \, :\, f(c)(b-a) \le \int_a^b{f(x) \, dx}\right\} \right)>0
\end{equation}
and
\begin{equation}
\label{pm2}
m\left(\left\{c \in (a,b) \, :\,  \int_a^b{f(x) \, dx} \le f(c)(b-a)\right\} \right)>0.
\end{equation}
\end{theorem}

\noindent
{\bf Proof.} The set
\begin{align*}
A&=\left\{c \in [a,b] \, :\, f(c)(b-a) \le \int_a^b{f(x) \, dx}\right\}\\
&=f^{-1}\left( \left[-\infty,\dfrac{1}{b-a}\int_a^b{f(x) \, dx}\right] \right)
\end{align*}
is Lebesgue measurable. If $m(A)=0$, i.e., if for almost all $c \in [a,b]$ we have
$$f(c) > \dfrac{1}{b-a}\int_a^b{f(x) \, dx},$$
then, integrating in both sides of this inequality, we obtain that
$$\int_a^b{f(x) \, dx}>\int_a^b{f(x) \, dx},$$
a contradiction. The proof of (\ref{pm2}) is similar. \qed

Theorem \ref{vmmvt} can be equivalently stated in terms of derivatives of absolutely continuous functions, and then it looks like one of those mean value theorems in differential calculus rather than in integral calculus. 
\begin{corollary}
\label{comvt}
If $F:I=[a,b]\longrightarrow \R$ is absolutely continuous then
\begin{equation}
\label{pm11}
m\left(\left\{c \in (a,b) \, :\, F'(c)(b-a) \le F(b)-F(a)\right\} \right)>0
\end{equation}
and
\begin{equation}
\label{pm21}
m\left(\left\{c \in (a,b) \, :\, F(b)-F(a) \le F'(c)(b-a)\right\} \right)>0.
\end{equation}
\end{corollary} 

\noindent
{\bf Proof.} Use Theorem \ref{vmmvt} with $f$ replaced by $F'$. \qed

\begin{remark}
In the conditions of Corollary \ref{comvt} the set of points $c$ such that
$$F(b)-F(a)=F'(c)(b-a)$$
may be empty. As an example, consider $F(x)=|x|$ for all $x \in [-1,1]$.

Analogously, in the conditions of Theorem \ref{vmmvt} there may be no point $c \in [a,b]$ satisfying
$$\int_a^b{f(x) \, dx}=f(c)(b-a).$$
\end{remark}

Part (C) in Theorem \ref{tkol} is now immediate from Corollary \ref{comvt}.

\begin{proposition}[Part (C) in Theorem \ref{tkol}]
If $F:I=[a,b]\longrightarrow \R$ is absolutely continuous and $F'=0$ almost everywhere in $I$ then $F$ is constant.
\end{proposition}

\noindent
{\bf Proof.} For each $x \in (a,b]$ we have, by (\ref{pm11}), 
$$m\left(\left\{c \in (a,x) \, :\, F'(c)(x-a) \le F(x)-F(a)\right\} \right)>0.$$
Since $F'=0$ almost everywhere, we conclude that $0 \le F(x)-F(a)$.

Similarly, we deduce from (\ref{pm21}) that $F(x)-F(a) \le 0$, and therefore we have $F(x)=F(a)$ for all $ x\in I$.
\qed

\section{Mean Value Theorems for Riemann integrals}
Our Mean Value Theorem \ref{vmmvt} yields the following Mean Value Theorem for the Riemann integral.
\begin{theorem}
\label{vmmvtr}
If $f:I=[a,b]\longrightarrow \R$ is Riemann--integrable then there exist points $c_i \in (a,b)$ ($i=1,2$) such that $f$ is continuous at $c_i$ ($i=1,2$) and
$$f(c_1)(b-a) \le \int_a^b{f(x) \, dx}\le f(c_2) (b-a).$$
\end{theorem}

\noindent
{\bf Proof.} It follows from (\ref{pm1}) that the set
$$\left\{c \in (a,b) \, :\, f(c)(b-a) \le \int_a^b{f(x) \, dx}\right\}$$
cannot be a part of the set of discontinuity points of $f$, because the latter is null. Hence there exists at least one $c_1$ in the conditions of the statement.

The existence $c_2$ follows from (\ref{pm2}). \qed

Theorem \ref{t1} is now a consequence of Theorem \ref{vmmvtr}.

\bigbreak

\noindent
{\bf Second proof of Theorem \ref{t1}}. Following the first proof of Theorem \ref{t1}, we know that $f=0$ whenever $f$ is continuous. Now let $x \in (a,b]$ be fixed and use Theorem \ref{vmmvtr} on the interval $[a,x]$ to deduce that
$$0 \le F(x)-F(a)=\int_a^x{f(y) \, dy}\le 0.$$
\qed

The rest of this section is devoted to proving Theorem \ref{t1} using only basic elements of Riemann integration. In particular, we even avoid using null--measure sets. We are convinced that the following material could be interesting for a broad part of the mathematical community, even for undergraduate students. The fundamentals of Riemann integration can be looked up in any textbook, and \cite{tre} is specially accessible and accurate.

\subsection{Riemann integrability yields some continuity}
Integrable functions are continuous at many points, and the usual way to prove it uses the notion of oscillation of a function in an interval and at a point. The following standard material is included for self--containedness.

\bigbreak

For a bounded function $f:[a,b]
\longrightarrow \R$ and a subinterval $[c,d] \subset [a,b]$ we call the oscillation of $f$ in the subinterval $[c,d]$ the number
$$\mbox{osc}(f,[c,d])=\sup_{c \le x \le d}f(x)-\inf_{c \le x \le d}f(x).$$
The oscillation fulfills the following three basic properties: $\mbox{osc}(f,[c,d]) �\ge 0$; for all $x,y \in [c,d]$ we have $|f(x)-f(y)| \le \mbox{osc}(f,[c,d])$;
 and if $[\hat c, \hat d] \subset [c,d]$ then $\mbox{osc}(f,[\hat c, \hat d]) \le \mbox{osc}(f,[c,d]).$

Subsequently, we define the oscillation of $f$ at a point $c \in (a,b)$ as
$$\mbox{osc}(f,c)=\lim_{\delta \to 0^+}\mbox{osc}(f,[c-\delta,c+\delta])=\inf_{\delta >0}\mbox{osc}(f,[c-\delta,c+\delta]),$$
and it is an exercise to prove that $f$ is continuous at $c \in (a,b)$ if $\mbox{osc}(f,c)=0$ (the converse is also true and easy to prove, but it is not essential for this paper).

\bigbreak

The following lemma on continuity of integrable functions at some points is somewhat na\"{\i}ve in comparison with the Lebesgue's test for Riemann integrability. The main reason for not invoking Lebesgue's test in this section is that we only need a very simple connection between Riemann integrability and continuity (which wants a simple proof, even adequate for an elementary course).

\begin{lemma}
\label{lemcon}
If $f:[a,b]\longrightarrow \R$ is Riemann--integrable on $[a,b]$ then there exists at least one point $c \in (a,b)$ at which $f$ is continuous.
\end{lemma}

\noindent
{\bf Proof.} Since $f$ is integrable on $[a,b]$ for every $\varepsilon >0$ there exists a partition $P=\{x_0,x_1,\dots,x_n\}$ such that
$$\varepsilon \, (b-a) >U(f,P)-L(f,P)=\sum_{k=1}^{n}{\mbox{osc}(f,[x_{k-1},x_k])}(x_k-x_{k-1}),$$
where $U(f,P)$ is the upper sum of $f$ relative to $P$ and $L(f,P)$ is the lower sum. Hence for some $j \in \{1,2,\dots,n\}$ we have $\mbox{osc}(f,[x_{j-1},x_j])<\varepsilon.$

Replacing $[x_{j-1},x_j]$ by one of its subintervals, if neccesary, we also have
$$[x_{j-1},x_j] \subset (a,b) \quad \mbox{and} \quad x_j-x_{j-1}<\varepsilon.$$

 Since $f$ is integrable on every subinterval of $[a,b]$ we can construct inductively a sequence $\{[a_n,b_n]\}_{n \in {\scriptsize \N}}$ of nested subintervals of $[a,b]$ such that for all $n \in \N$, $n \ge 2$, we have
$$\mbox{osc}(f,[a_n,b_n])<1/n, \, \, \, \, [a_n,b_n] \subset (a_{n-1},b_{n-1}), \, \, \, \,  \mbox{and} \quad b_n-a_n<1/n.$$ 

Therefore there is a unique point $c \in \cap_{n \in {\scriptsize \N}}[a_n,b_n] \subset (a,b)$, and for all $n \in \N$ we have
$$0 \le \mbox{osc}(f,c)�\le \mbox{osc}(f,[a_n,b_n])<\dfrac{1}{n},$$
which implies that $f$ is continuous at $c$.
\qed

Our next corollary reveals that Lemma \ref{lemcon} is sharper than it might look.

\begin{corollary}
\label{corlem}
If $f:[a,b]\longrightarrow \R$ is Riemann--integrable on $[a,b]$ then $f$ is continuous in a dense subset of $[a,b]$.
\end{corollary}

\noindent
{\bf Proof.} Use Lemma \ref{lemcon} in each nondegenerate subinterval $[x,y] \subset [a,b]$. \qed

\subsection{Proving Theorem \ref{t1} without null--measure sets}
 We base our third and last proof on the following lemma, which is the closest we can get to Theorem \ref{vmmvtr} without using null--measure sets.

\begin{lemma}
\label{lem1}
If $f:[a,b] \longrightarrow \R$ is Riemann--integrable on $[a,b]$ then for every $\varepsilon>0$ there exist points $c_i \in (a,b)$ ($i=1,2$) such that $f$ is continuous at $c_i$ ($i=1,2$),
\begin{equation}
\label{fi}
f(c_1)(b-a)<\int_a^b{f(x) \, dx}+\varepsilon
\end{equation}
and
\begin{equation}
\label{fi2}
\int_a^b{f(x) \, dx}-\varepsilon< f(c_2)(b-a).
\end{equation}
\end{lemma}

\noindent
{\bf Proof.} Let $\varepsilon>0$ be fixed and let $P=\{x_0,x_1,\dots,x_n\}$ be a partition of $[a,b]$ such that
$$ U(f,P)=\sum_{k=1}^n{\sup_{x_{k-1}\le x \le x_k}f(x) (x_k-x_{k-1})} <\int_a^b{f(x) \, dx}+ \varepsilon.$$
This implies the existence of some $j \in \{1,2,\dots,n\}$ such that
$$\sup_{x_{j-1}\le x\le x_j}f(x) < \dfrac{1}{b-a} \left(\int_a^b{f(x) \, dx}+ \varepsilon \right).$$
Now (\ref{fi}) is satisfied with any $c_1 \in (x_{j-1},x_j)$ such that $f$ is continuous at $c_1$ (such $c_1$ exists by virtue of Corollary \ref{corlem}).

The proof of (\ref{fi2}) is similar and involves lower sums.\qed

\noindent
{\bf Third proof to Theorem \ref{t1}}. Following our first proof of Theorem \ref{t1} we know that $f=0$ whenever $f$ is continuous. Now for each $x \in (a,b]$ and $\varepsilon>0$ we use Lemma \ref{lem1} to guarantee the existence of some $c_1 \in (a,x)$ such that $f$ is continuous at $c_1$ and
$$F(x)-F(a)+\varepsilon=\int_a^x{f(y) \, dy}+\varepsilon >f(c_1)(x-a)=0.$$
Since $\varepsilon>0$ was arbitrary, we conclude that $F(x)-F(a) \ge 0$. The reverse inequality is deduced from Lemma \ref{lem1} in a similar way.
\qed

\subsection{Some more mean value theorems}
This section collects the most basic Mean Value Theorems in this paper, and we are not going to use them in connection with Theorem \ref{t1}. Despite some of them are known, they are not easily traceable in the literature and that is why we have decided to include them here.

\bigbreak

It follows from the definitions that any bounded function $f:[a,b] \longrightarrow \R$ satisfies
\begin{equation}
\label{mvin}
m(b-a) \le \underline{\int_a^b}{f(x) \, dx} \le \overline{\int_a^b}{f(x) \, dx}\le M(b-a),
\end{equation}
where $m=\inf \{f(x) \, : \, x \in [a,b]\}$ and $M=\sup \{f(x) \, : \, x \in [a,b]\}.$

The well--known Mean Value Theorems for the Riemann integral are immediate consequences of (\ref{mvin}), but a deeper analysis leads to better information. We start proving the following mean value inequalities, which improve (\ref{mvin}).
\begin{theorem}
\label{mvsin}
Let $f:[a,b] \longrightarrow \R$ be bounded in $[a,b]$. 

If $f$ is continuous at one point in $(a,b)$ then there exist points $c_1, \, c_2 \in (a,b)$ such that
\begin{equation}
\label{mvin2}
f(c_1)(b-a) \le \underline{\int_a^b}{f(x) \, dx} \le \overline{\int_a^b}{f(x) \, dx}\le f(c_2)(b-a).
\end{equation}
\end{theorem}

\noindent
{\bf Proof.} We shall only prove that there is some $c_2 \in (a,b)$ satisfying the right--hand inequality in (\ref{mvin2}), as the proof is analogous for the left--hand inequality.

Reasoning by contradiction, we assume that for all $x \in (a,b)$ we have
\begin{equation}
\label{in}
f(x)<\dfrac{1}{b-a}\overline{\int_a^b}{f(y) \, dy}.
\end{equation}
Changing the values of $f$ at $a$ and $b$, if necessary (which does not alter the value of the upper integral), we have a new function $\hat f$ such that (\ref{in}) holds for all $x \in [a,b]$ with $f$ replaced by $\hat f$.

Since $\hat f$ is continuous at some $c \in (a,b)$, we can find $\varepsilon>0$ so that
$$\hat f(x)<\dfrac{1}{b-a}\overline{\int_a^b}{\hat f(y) \, dy}-\varepsilon \quad \mbox{for $x \in (c-\varepsilon,c+\varepsilon) \subset (a,b).$}$$
Finally, consider the partition $P=\{a,c-\varepsilon,c+\varepsilon,b\}$ to obtain
$$\overline{\int_a^b}{\hat f(x) \, dx}\le U(\hat f,P) \le \overline{\int_a^b}{\hat f(x) \, dx}-2 \varepsilon^2< \overline{\int_a^b}{\hat f(x) \, dx},$$
a contradiction. \qed

\begin{remark}
The continuity condition cannot be omitted in Theorem \ref{mvsin}. 

Indeed, let $\{p_n/q_n\}_{n�\in {\scriptsize \N}}$ be a sequential arrangement of all rational numbers in $[0,1]$, with $p_n, q_n$ positive integers and $p_n/q_n$ irreductible for all $n \in \N$. The function $f:[0,1] \longrightarrow \R$ defined by
$$f(p_n/q_n)=1/q_n, \quad f(\pi \, p_n/q_n)=1-1/q_n, \quad \mbox{and $f(x)=1/2$ otherwise,}$$
satisfies $0<f(x)<1$ for all $x \in [0,1]$,
$$\underline{\int_0^1}{f(x) \, dx}=0 \quad \mbox{and} \quad \overline{\int_a^b}{f(x) \, dx}=1.$$
\end{remark}

Integrability ensures continuity to some extent, and therefore the parti\-cu\-lar case of Theorem \ref{mvsin} for integrable functions has a cleaner statement.

\begin{corollary}
\label{corpro}
If $f:[a,b] \longrightarrow \R$ is Riemann--integrable on $[a,b]$ then there exist points $c_1,\, c_2 \in (a,b)$ such that
\begin{equation}
\label{mvin3}
f(c_1)(b-a) \le  {\int_a^b}{f(x) \, dx}  \le f(c_2)(b-a).
\end{equation}

In particular, if $f$ is continuous on $[a,b]$ then there exists at least one $c \in (a,b)$ such that
\begin{equation}
\label{mvincont}
\int_a^b{f(x) \, dx}=f(c)(b-a).
\end{equation}

\end{corollary}

\noindent
{\bf Proof.} Inequality (\ref{mvin3}) follows from Theorem \ref{mvsin} and Lemma \ref{lemcon}. Using Bolzano's Theorem, we deduce (\ref{mvincont}) from (\ref{mvin3}). \qed

\begin{remark}
As far as the author knows, the first part of Corollary \ref{corpro}, concerning integrable functions, is new. Notice that it has some interesting consequences, such as the strict positivity of the integral of strictly positive integrable functions (which, in turn, yields the strict monotonicity of the Riemann integral).

The Mean Value Theorem for continuous functions in Corollary \ref{corpro} is known, but difficult to find in textbooks. We can cite Stromberg \cite{str}, where we find it on page 281, left as Exercise 28. The fact that (\ref{mvincont}) is satisfied with some $c$ in the open interval $(a,b)$ is interesting and useful. For instance, as observed by \'Oscar L\'opez Pouso, it allows us to pass from the integral expression of Taylor's remainder to its usual differential form.
\end{remark}

\end{document}